\documentclass{article}
\usepackage[english]{babel}
\usepackage{latexsym,amssymb,amsmath}
\usepackage[cp1251]{inputenc}
\usepackage{amsfonts,amssymb}
\usepackage{graphicx,graphics,hhline}
\usepackage{euscript}

\begin{document}

\begin{center}
{\Large \textbf{On $h(x)$ -- Fibonacci polynomials in an arbitrary algebra}}

\begin{equation*}
\end{equation*}%
\textbf{Cristina Flaut, Vitalii Shpakivskyi and Elena Vlad }
\end{center}

\begin{equation*}
\end{equation*}

\textbf{Abstract. }{\small In this paper, we introduce }$h(x)${\small \ --
Fibonacci polynomials in an arbitrary finite-dimensional unitary algebra
over a field }${\small K\ }\left( K=\mathbb{R},\mathbb{C}\right) ,${\small \
which generalize both }$h(x)${\small \ -- Fibonacci quaternion polynomials
and }$h(x)${\small \ -- Fibonacci octonion polynomials. For }$h(x)${\small \
-- Fibonacci polynomials in such an arbitrary algebra, we prove summation
formula, generating function, Binet-style formula, Catalan-style identity,
and d'Ocagne-type identity. }

\begin{equation*}
\end{equation*}

MSC 2010: Primary 11B39; Secondary 11R54%
\begin{equation*}
\end{equation*}

\textbf{Introduction}

\bigskip

\begin{equation*}
\end{equation*}%
In modern investigation there is a huge interest to real Fibonacci numbers
and numerous their generalizations. The Fibonacci numbers $f_{n}$ are the
terms of the sequence $0,1,1,2,3,5,\ldots $, where 
\begin{equation*}
f_{n}=f_{n-1}+f_{n-2},\quad n=2,3,\ldots ,
\end{equation*}%
with the initial values $f_{0}=0,f_{1}=1$.

In the paper \cite{Falcon}, were introduced $k$--Fibonacci numbers $f_{k,n}$
by the equality 
\begin{equation*}
f_{k,n}=kf_{k,n-1}+f_{k,n-2},\quad n=2,3,\ldots ,
\end{equation*}%
with the initial values $f_{k,0}=0,f_{k,1}=1,$ for all $k$. It is easy to
see that Pell numbers are $2$--Fibonacci numbers. Catalan studied
polynomials $F_{n}(x),$ defined by the recurrence relation 
\begin{equation*}
F_{n}(x)=xF_{n-1}(x)+F_{n-2}(x),\quad n\geq 3,
\end{equation*}%
where $F_{1}(x)=0$, $F_{2}(x)=x$. There are another generalizations of
Fibonacci numbers, as for example, Jacobsthal polynomials $J_{n}(x)$,
defined by 
\begin{equation*}
J_{n}(x)=J_{n-1}(x)+xJ_{n-2}(x),\quad n=3,4,\ldots ,
\end{equation*}%
where $J_{1}(x)=J_{2}(x)=1$. There also exist \ the Byrd polynomials $%
\varphi _{n}(x)$ which are defined by the relation 
\begin{equation*}
\varphi _{n}(x)=2x\varphi _{n-1}(x)+\varphi _{n-2}(x),\quad n=2,3,\ldots ,
\end{equation*}%
where $J_{0}(x)=0,\varphi _{1}(x)=1$. The Lucas polynomials $L_{n}(x)$ were
introduction by Bicknell, and are defined by 
\begin{equation*}
L_{n}(x)=xL_{n-1}(x)+L_{n-2}(x),\quad n=2,3,\ldots ,
\end{equation*}%
where $L_{0}(x)=2,L_{1}(x)=x$.

In the paper \cite{Nalli}, authors introduced $h(x)$--Fibonacci polynomials $%
F_{h,n}(x),$ which generalize both Catalan's polynomials, Byrd's polynomials
and also the $k$-Fibonacci numbers. Let $h(x)$ be a polynomial with real
coefficients. The $h(x)$--Fibonacci polynomials are defined by the
recurrence relation 
\begin{equation}
F_{h,n}(x)=h(x)F_{h,n-1}(x)+F_{h,n-2}(x),\quad n=2,3,\ldots ,  \label{1}
\end{equation}%
where $F_{h,0}(x)=0$, $F_{h,1}(x)=1,$ for all polynomials $h(x)$. In \cite%
{Nalli}, are studied basic properties of $h(x)$--Fibonacci polynomials. We
note that due to essential results given in the paper \cite{Nalli} were
obtained the results from the papers \cite{Falcon-07,Falcon-09}, where were
considered the $x$--Fibonacci polynomials.

Numerous generalizations of Fibonacci numbers generated different
hypercomplex generalizations of Fibonacci numbers. Therefore, Fibonacci
elements over some special algebras were intensively studied in the last
time in various papers, as for example: \cite{Akk}~---~\cite{Sw-73}. All
these papers studied properties of Fibonacci elements in complex numbers, or
in quaternions and octonions, or in generalized Quaternion and Octonion
algebras, or studied dual vectors or dual Fibonacci quaternions. At the same
time, in the paper \cite{Flaut-Shpa-Bull} were considered Fibonacci elements
in an arbitrary finite-dimensional unitary algebra over a field $K$ ($K=%
\mathbb{R},\mathbb{C}$) and were proved some basic properties of these
hypercomplex numbers (generating functions, Binet formula, Cassini's
identity, etc.).

In the paper \cite{Ra-15} are defined the $k$ -- Fibonacci and the $k$ --
Lucas quaternions. For these quaternions were investigated the generating
functions, Binet formula, formulae for some sums and Cassini's identity.
Some of results of the paper \cite{Ra-15} were generalized in the article 
\cite{Catarino}, where was introduced the $h(x)$ -- Fibonacci quaternion
polynomials which generalize the $k$ -- Fibonacci quaternion numbers. In 
\cite{Catarino}, is presented a Binet-style formula, ordinary generating
function and some basic identities for $h(x)$ -- Fibonacci quaternion
polynomials.

In the paper \cite{Ipec} are defined $h(x)$ -- Fibonacci octonion
polynomials. For the last mentioned, were obtained a similar Binet formula
and generating function.

In this paper, we introduce $h(x)$ -- Fibonacci polynomials in an arbitrary
finite dimensional unitary algebra over a field $K$ ($K=\mathbb{R},\mathbb{C}
$) which generalize both $k$ -- Fibonacci quaternions, $h(x)$ -- Fibonacci
quaternion polynomials and $h(x)$ -- Fibonacci octonion polynomials. We also
prove some relations between $h(x)$ -- Fibonacci polynomials in such an
arbitrary algebra, Binet-style formula, Catalan-style identity, and
generating function. 
\begin{equation*}
\end{equation*}%
\textbf{Real }$h(x)$\textbf{\ -- Fibonacci polynomials and their properties}%
\begin{equation*}
\end{equation*}

In this section we indicate some basic properties of $h(x)$ -- Fibonacci
polynomials defined by the equality (\ref{1}).

1. Generating function \cite{Nalli}:\quad 
\begin{equation*}
g(t)=\frac{t}{1-h(x)t-t^{2}}.
\end{equation*}

2. For $n\in N$ (see \cite{Nalli}), we have 
\begin{equation}
F_{h,n}(x)=\sum_{k=0}^{[\frac{n-1}{2}]}\complement
_{n-k-1}^{k}\,h^{n-2k-1}(x).  \label{2.0}
\end{equation}

3. For $n\in N$ (see \cite{Nalli}), we have 
\begin{equation*}
F_{h,n}(x)=2^{1-n}\sum_{k=0}^{[\frac{n-1}{2}]}\complement
_{n}^{2k+1}\,h^{n-2k-1}(x)(h^{2}(x)+4)^{k}.
\end{equation*}

4. For $n\in N$ (see \cite{Nalli}), we have 
\begin{equation*}
F_{h,n}(x)=i^{n-1}U_{n-1}\left( \frac{h(x)}{2i}\right) ,
\end{equation*}%
where $i^{2}=-1$ and $U_{n}(t):=\sum\limits_{j=0}^{[\frac{n-1}{2}%
]}(-1)^{j}\,\complement _{n-j}^{j}(2t)^{n-2j}$ is the Chebyshev polynomial
of the second kind.

5. Let $\alpha (x)$ and $\beta (x)$ denote the roots of the characteristic
equation $v^{2}-h(x)v-1=0$ of the recurrence relation (\ref{1}). Then 
\begin{equation}
\alpha (x):=\frac{h(x)+\sqrt{h^{2}(x)+4}}{2}\,,\quad \beta (x):=\frac{h(x)-%
\sqrt{h^{2}(x)+4}}{2}\,.  \label{2}
\end{equation}%
For all $n=0,1,2,\ldots $ the Binet formula is of the form (see \cite{Nalli}%
) 
\begin{equation}
F_{h,n}(x)=\frac{\alpha ^{n}(x)-\beta ^{n}(x)}{\alpha (x)-\beta (x)}\,.
\label{Binet}
\end{equation}

6. (see \cite{Falcon-09}): \quad $\displaystyle\lim\limits_{n\rightarrow
\infty }\frac{F_{h,n+1}(x)}{F_{h,n}(x)}=\alpha (x)$.

7. (see \cite{Falcon-09}): \quad 
\begin{equation}
\sum_{k=1}^{n}F_{h,k}(x)=\frac{F_{h,n+1}(x)+F_{h,n}(x)-1}{h(x)}\,.  \label{3}
\end{equation}

8. For $n,r$ integers and $n>r$ we have the Catalan identity \cite{Falcon-09}%
: 
\begin{equation*}
F_{h,n-r}(x)F_{h,n+r}(x)-F_{h,n}^{2}(x)=(-1)^{n-r-1}F_{h,r}^{2}(x).
\end{equation*}

9. For $a,b,c,d$ and $r$ integers, with $a+b=c+d$ we have \cite{Falcon-09}: 
\begin{equation*}
F_{h,a}(x)F_{h,b}(x)-F_{h,c}(x)F_{h,d}(x)=(-1)^{r}\big(%
F_{h,a-r}(x)F_{h,b-r}(x)-F_{h,c-r}(x)F_{h,d-r}(x)\big).
\end{equation*}

We note that the representation given in (\ref{2.0}) can be rewritten in the
following differential form

10. 
\begin{equation*}
F_{h,n}(x)=\sum_{k=0}^{[\frac{n-1}{2}]}\frac{1}{k!}\,\frac{d^{k}}{dh^{k}}%
\,h^{n-k-1}(x).
\end{equation*}
\begin{equation*}
\end{equation*}

$h(x)$\textbf{\ -- Fibonacci polynomials in an arbitrary finite dimensional
algebra}%
\begin{equation*}
\end{equation*}

Let $A$ be an unitary arbitrary $(m+1)$-dimensional algebra over $K$ ($K=%
\mathbb{R},\mathbb{C}$) with a basis $\{e_{0},e_{1},e_{2},...,e_{m}\}.%
\medskip $

\textbf{Definition 3.1. }The $h(x)$ -- Fibonacci polynomials $%
\{Q_{h,n}(x)\}_{n=0}^{\infty }$ in such an arbitrary algebra $A$ are defined
by the recurrent relation 
\begin{equation}
Q_{h,n}(x)=\sum\limits_{k=0}^{m}F_{h,n+k}(x)\,e_{k},  \label{ozn-h-Alg}
\end{equation}%
where $F_{h,n}(x)$ is the $n$th real $h(x)$ -- Fibonacci polynomial.

In the case where the algebra $A$ coincides with the quaternion algebra $%
\mathbb{H}$, we obtain $h(x)$ -- Fibonacci quaternion polynomials, studied
in \cite{Catarino}. If an algebra $A$ coincides with the octonion algebra $H$%
, we obtain $h(x)$ -- Fibonacci octonion polynomials which were considered
in \cite{Ipec}.\vskip2mm

\textbf{Proposition 3.2.} \textit{For any natural numbers} $n$ \textit{and} $%
p$ \textit{the following relations hold:}\vskip2mm

$(i)$ 
\begin{equation*}
Q_{h,n+2}(x)=h(x)Q_{h,n+1}(x)+Q_{h,n}(x);
\end{equation*}

$(ii)$ 
\begin{equation*}
\overset{p}{\underset{k=1}{\sum }}Q_{h,k}(x)=\frac{1}{h(x)}%
(Q_{h,p+1}(x)+Q_{h,p}(x)-Q_{h,0}(x)-Q_{h,1}(x)).
\end{equation*}

\textbf{Proof.} $(i)$ directly follows from definition 3.1. $(ii)$ In the
following, instead of $Q_{h,n}(x)$ and $F_{h,j}(x)$ we will write $Q_{h,n}$
and $F_{h,j}$, respectively. Using the identity (\ref{3}), we have 
\begin{equation*}
\sum\limits_{k=1}^{p}Q_{h,k}=\sum\limits_{j=0}^{m}F_{h,j+1}\,e_{j}+\sum%
\limits_{j=0}^{m}F_{h,j+2}\,e_{j}+\cdots
+\sum\limits_{j=0}^{m}F_{h,j+p}\,e_{j}=
\end{equation*}%
\begin{equation*}
e_{0}(F_{h,1}+F_{h,2}+\cdots +F_{h,p})+e_{1}(F_{h,2}+F_{h,3}+\cdots
+F_{h,p+1})+\cdots
\end{equation*}%
\begin{equation*}
+e_{m}(F_{h,m+1}+F_{h,m+2}+\cdots +F_{h,m+p})=\frac{e_{0}}{h(x)}%
(F_{h,p+1}+F_{h,p}-1)+
\end{equation*}%
\begin{equation*}
\frac{e_{1}}{h(x)}(F_{h,p+2}+F_{h,p+1}-1-h(x)F_{h,1})+\cdots
\end{equation*}%
\begin{equation*}
+\frac{e_{m}}{h(x)}(F_{h,p+m+1}+F_{h,p+m}-1-h(x)F_{h,1}-h(x)F_{h,2}-\cdots
-h(x)F_{h,m}).
\end{equation*}

Since from (\ref{3})\thinspace\ we have $1+h(x)\sum%
\limits_{k=1}^{r}F_{h,k}=F_{h,r}+F_{h,r+1}$, then we get 
\begin{equation*}
\sum\limits_{k=1}^{p}Q_{h,k}=\frac{e_{0}}{h(x)}(F_{h,p+1}+F_{h,p}-1)+\frac{%
e_{1}}{h(x)}(F_{h,p+2}+F_{h,p+1}-F_{h,1}-F_{h,2})+\cdots
\end{equation*}%
\begin{equation*}
+\frac{e_{m}}{h(x)}(F_{h,p+m+1}+F_{h,p+m}-F_{h,m}-F_{h,m+1})=\frac{1}{h(x)}%
(Q_{h,p+1}+Q_{h,p}-Q_{h,0}-Q_{h,1}).
\end{equation*}%
\noindent The proposition is proved. \vskip2mm

We obtain the following Binet formula for $Q_{h,n}(x)$.\vskip2mm

\textbf{Theorem 3.3.} \textit{For} $n=0,1,2,\ldots $ \textit{we have the
following relation} 
\begin{equation}
Q_{h,n}(x)=\frac{\alpha ^{\ast }(x)\alpha ^{n}(x)-\beta ^{\ast }(x)\beta
^{n}(x)}{\alpha (x)-\beta (x)},  \label{Binet-hyp}
\end{equation}%
\textit{where} $\alpha ^{\ast }=\sum\limits_{k=0}^{m}\alpha
^{k}(x)\,e_{k}\,,\quad \beta ^{\ast }=\sum\limits_{k=0}^{m}\beta
^{k}(x)\,e_{k}\,.$\vskip2mm

\textbf{Proof.} Using the Binet-style formula (\ref{Binet}), we obtain 
\begin{equation*}
Q_{h,n}(x)=\sum\limits_{k=0}^{m}F_{h,n+k}(x)\,e_{k}=\frac{\alpha
^{n}(x)-\beta ^{n}(x)}{\alpha (x)-\beta (x)}\,e_{0}+\frac{\alpha
^{n+1}(x)-\beta ^{n+1}(x)}{\alpha (x)-\beta (x)}\,e_{1}+\cdots
\end{equation*}%
\begin{equation*}
+\frac{\alpha ^{n+m}(x)-\beta ^{n+m}(x)}{\alpha (x)-\beta (x)}\,e_{m}=\frac{%
\alpha ^{n}(x)}{\alpha (x)-\beta (x)}(e_{0}+\alpha (x)e_{1}+\alpha
^{2}(x)e_{2}+\cdots
\end{equation*}%
\begin{equation*}
+\alpha ^{m}(x)e_{m})-\frac{\beta ^{n}(x)}{\alpha (x)-\beta (x)}\left(
e_{0}+\beta (x)e_{1}+\beta ^{2}(x)e_{2}+\cdots +\beta ^{m}(x)e_{m}\right) =
\end{equation*}%
\begin{equation*}
\frac{\alpha ^{\ast }(x)\alpha ^{n}(x)-\beta ^{\ast }(x)\beta ^{n}(x)}{%
\alpha (x)-\beta (x)}\,.
\end{equation*}%
\vskip2mm

\textbf{Remark 3.4.} The above result generalizes the Binet formulae from
the papers \cite{Catarino} and \cite{Ipec}.\vskip2mm

\textbf{Definition 3.5.} The generating function $G(t)$ of the sequence $%
\{Q_{h,n}(x)\}_{n=0}^{\infty }$ is defined by 
\begin{equation}
G(t)=\sum\limits_{n=0}^{\infty }Q_{h,n}(x)t^{n}.  \label{gen-fun}
\end{equation}%
\vskip2mm

\textbf{Theorem 3.6.} \textit{The generating function for the} $h(x)$ 
\textit{-- Fibonacci polynomials} $Q_{h,n}(x)$ \textit{in an arbitrary
algebra is of the form} 
\begin{equation*}
G(t)=\frac{Q_{h,0}(x)+(Q_{h,1}(x)-h(x)Q_{h,0}(x))t}{1-h(x)t-t^{2}}.
\end{equation*}%
\vskip2mm

\textbf{Proof. }Taking into account the equality (\ref{gen-fun}), we
consider the product 
\begin{equation*}
G(t)(1-h(x)t-t^{2})=\sum\limits_{n=0}^{\infty
}Q_{h,n}(x)t^{n}-h(x)\sum\limits_{n=0}^{\infty
}Q_{h,n}(x)t^{n+1}-\sum\limits_{n=0}^{\infty }Q_{h,n}(x)t^{n+2}=
\end{equation*}%
\begin{equation*}
Q_{h,0}(x)+(Q_{h,1}(x)-h(x)Q_{h,0}(x)+\sum\limits_{n=2}^{\infty
}t^{n}(Q_{h,n}-h(x)Q_{h,n-1}-Q_{h,n-2})=
\end{equation*}%
\begin{equation*}
Q_{h,0}(x)+(Q_{h,1}(x)-h(x)Q_{h,0}(x)).
\end{equation*}%
The theorem is proved. \vskip2mm

\textbf{Remark 3.7.} The above Theorem generalizes results from the papers 
\cite{Catarino} and \cite{Ipec}.\vskip2mm

\textbf{Theorem 3.8.} (Catalan's identity) \textit{For nonnegative integer
numbers} $n,r$, \textit{such that} $r\leq n$, \textit{we have} 
\begin{equation*}
Q_{h,n+r}(x)Q_{h,n-r}(x)-Q_{h,n}^{2}(x)=
\end{equation*}%
\begin{equation*}
\frac{(-1)^{n+r+1}}{h^{2}(x)+4}(\alpha ^{\ast }(x)\beta ^{\ast
}(x)[(-1)^{r+1}+\alpha ^{2}(x)]+\beta ^{\ast }(x)\alpha ^{\ast
}(x)[(-1)^{r+1}+\beta ^{2}(x)]).
\end{equation*}

\vskip2mm \textbf{Proof.} Using the formula (\ref{Binet-hyp}), we have 
\begin{equation*}
Q_{h,n+r}(x)Q_{h,n-r}(x)-Q_{h,n}^{2}(x)=
\end{equation*}%
\begin{equation*}
\frac{1}{(\alpha (x)-\beta (x))^{2}}(\alpha ^{\ast }(x)\beta ^{\ast
}(x)(\alpha (x)\beta (x))^{n}\left[ 1-\left( \frac{\alpha (x)}{\beta (x)}%
\right) ^{r}\,\right] +
\end{equation*}%
\begin{equation*}
\beta ^{\ast }(x)\alpha ^{\ast }(x)(\alpha (x)\beta (x))^{n}\left[ 1-\left( 
\frac{\beta (x)}{\alpha (x)}\right) ^{r}\,\right] ).
\end{equation*}%
Now, taking into account the relations $\alpha (x)\beta (x)=-1$, $\frac{%
\alpha (x)}{\beta (x)}=-\alpha ^{2}(x)$, $\frac{\beta (x)}{\alpha (x)}%
=-\beta ^{2}(x)$, we obtain the statement of the theorem. The theorem is now
proved. \vskip2mm

If in the Theorem 3.8 we set $r=1$, we obtain the Cassini-style identity.%
\vskip2mm

\textbf{Corollary 3.9.} (Cassini's identity) \textit{For any natural number }%
$n$\textit{, we have} 
\begin{equation*}
Q_{h,n+1}(x)Q_{h,n-1}(x)-Q_{h,n}^{2}(x)=
\end{equation*}%
\begin{equation*}
\frac{(-1)^{n}}{h^{2}(x)+4}(\alpha ^{\ast }(x)\beta ^{\ast }(x)[1+\alpha
^{2}(x)]+\beta ^{\ast }(x)\alpha ^{\ast }(x)[1+\beta ^{2}(x)]).
\end{equation*}%
\vskip2mm

\textbf{Remark 3.10. } Theorem 3.8 and Corollary 3.9 generalize the Theorems
3.8 and 3.9, respectively, from the paper \cite{Catarino}.\vskip2mm

Similarly to proof of Theorem 3.8, can be proved the following result. \vskip%
2mm

\textbf{Theorem 3.11.} (d'Ocagne's identity) \textit{Suppose that} $n$ 
\textit{is a nonnegative integer number and }$r$ \textit{is any natural
number. Then for} $r>n$, \textit{we have} 
\begin{equation*}
Q_{h,r}(x)Q_{h,n+1}(x)-Q_{h,r+1}(x)Q_{h,n}(x)=
\end{equation*}%
\begin{equation*}
\frac{(-1)^{n}}{\alpha (x)-\beta (x)}(\alpha ^{\ast }(x)\beta ^{\ast
}(x)\alpha ^{r-n}(x)-\beta ^{\ast }(x)\alpha ^{\ast }(x)\beta ^{r-n}(x)).
\end{equation*}

\textbf{Acknowledgement.} The second author is partially supported by Grant
of Ministry of Education and Science of Ukraine (Project No. 0116U001528).

\vskip4mm 
Cristina Flaut\newline
Faculty of Mathematics and Computer Science,\newline
Ovidius University,\newline
Bd. Mamaia 124, 900527, Constanta,\newline
Romania\newline
http://cristinaflaut.wikispaces.com/ \newline
http://www.univ-ovidius.ro/math/

e-mail: cflaut@univ-ovidius.ro,\ cristina$_{-}$flaut@yahoo.com

\vskip4mm

Vitalii Shpakivskyi\newline
Department of Complex Analysis and Potential Theory,\newline
Institute of Mathematics of the National Academy of Sciences of Ukraine,%
\newline
3, Tereshchenkivs'ka st.\newline
01601 Kyiv-4,\newline
Ukraine\newline
http://www.imath.kiev.ua/

e-mail: shpakivskyi86@gmail.com,\ shpakivskyi@imath.kiev.ua

\vskip4mm

Elena Vlad\newline
National College "Emil Botta", Adjud, Vrancea, Romania

elena\_brinzoi2000@yahoo.com%
\begin{equation*}
\end{equation*}

\end{document}